\documentclass[aos,MSNbibl,seceqn,nameyear,dvips]{arximspdf}
\usepackage{graphicx}

\doi{10.1214/11-AOS898}
\volume{39}
\issue{5}
\pubyear{2011}
\firstpage{2330}
\lastpage{2355}

\makeatletter

\newtheorem{theorem}{Theorem}[section]
\newtheorem{proposition}[theorem]{Proposition}
\newproclaim{Remark}{Remark}

\newcommand{\bR}{\mathbf R}
\newcommand{\bE}{\mathbb E}
\newcommand{\Dcal}{\mathcal D}
\newcommand{\Hcal}{\mathcal H}
\newcommand{\Lcal}{\mathcal L}
\newcommand{\Mcal}{\mathcal M}
\newcommand{\Pcal}{\mathcal P}
\newcommand{\Tcal}{\mathcal T}
\newcommand{\Wcal}{\mathcal W}

\newcommand{\var}{\operatorname{{\mathbb Var}}}
\newcommand{\argmin}{\mathop{\arg\min}}

\makeatother

\begin{document}
\begin{frontmatter}

\title{Optimal estimation of the mean function based on discretely sampled
functional data: Phase~transition}
\runtitle{Mean function estimation}

\begin{aug}
\author[A]{\fnms{T. Tony} \snm{Cai}\thanksref{t1}\ead[label=e1]{tcai@wharton.upenn.edu}}
\and
\author[B]{\fnms{Ming} \snm{Yuan}\corref{}\thanksref{t2}\ead[label=e2]{myuan@isye.gatech.edu}}
\runauthor{T. T. Cai and M. Yuan}
\affiliation{University of Pennsylvania and Georgia Institute of Technology}
\address[A]{Department of Statistics\\
The Wharton School\\
University of Pennsylvania\\
Philadelphia, Pennsylvania 19104\\
USA\\
\printead{e1}}
\address[B]{School of Industrial and Systems Engineering\\
Georgia Institute of Technology\\
Atlanta, Georgia 30332\\
USA\\
\printead{e2}} 
\end{aug}

\thankstext{t1}{Supported in part by NSF FRG Grant
DMS-08-54973.}

\thankstext{t2}{Supported in part by NSF Career Award DMS-08-46234.}

\received{\smonth{5} \syear{2010}}
\revised{\smonth{5} \syear{2011}}

%
\begin{abstract}
The problem of estimating the mean of random functions based on
discretely sampled data arises naturally in functional data analysis.
In this paper, we study optimal estimation of the mean function under
both common and independent designs. Minimax rates of convergence are
established and easily implementable rate-optimal estimators are
introduced. The analysis reveals interesting and different phase
transition phenomena in the two cases. Under the common design, the
sampling frequency solely determines the optimal rate of convergence
when it is relatively small and the sampling frequency has no effect on
the optimal rate when it is large. On the other hand, under the
independent design, the optimal rate of convergence is determined
jointly by the sampling frequency and the number of curves when the
sampling frequency is relatively small. When it is large, the sampling
frequency has no effect on the optimal rate. Another interesting
contrast between the two settings is that smoothing is necessary under
the independent design, while, somewhat surprisingly, it is not
essential under the common design.
\end{abstract}

%
\begin{keyword}
\kwd{Functional data}
\kwd{mean function}
\kwd{minimax}
\kwd{rate of convergence}
\kwd{phase transition}
\kwd{reproducing kernel Hilbert space}
\kwd{smoothing splines}
\kwd{Sobolev space}.
\end{keyword}

\end{frontmatter}

\section{Introduction}

Estimating the mean function based on discretely sampled
noisy observations is one of the most basic problems in functional
data analysis. Much progress has been made on developing estimation
methodologies. The two monographs by Ramsay and
Silverman (\citeyear{rs02}, \citeyear{rs05}) provide comprehensive
discussions on the
methods and applications. See also \citet{fv06}.

Let $X(\cdot)$ be a random function defined on the unit interval
$\Tcal=[0,1]$ and $X_1,\ldots,X_n$ be a sample of $n$ independent
copies of $X$. The goal is to estimate the mean function
$g_0(\cdot):=\bE(X(\cdot))$ based on noisy observations from discrete
locations on these curves:
%
%
\begin{equation}
\label{eqobsmod}
Y_{ij}=X_i(T_{ij})+\varepsilon_{ij},\qquad j=1,2,\ldots, m_i
\mbox{ and } i=1,2,\ldots, n,
\end{equation}
where $T_{ij}$ are sampling points, and $\varepsilon_{ij}$ are
independent random noise variables with $\bE\varepsilon_{ij}=0$ and
finite second moment $\bE\varepsilon_{ij}^2=\sigma_0^2<+\infty$. The sample
path of $X$ is assumed to be smooth in that it belongs to the usual
Sobolev--Hilbert spaces of order $r$ almost surely, such that
%
%
\begin{equation}
\label{eqsamplesmooth}
\bE\biggl( \int_\Tcal\bigl[X^{(r)}(t)\bigr]^2\,dt\biggr)<+\infty.
\end{equation}
Such problems naturally arise in a variety of applications and are
typical in functional data analysis [see, e.g., \citet
{rs05}, \citet{fv06}]. Various methods have been
proposed. However, little is known about their theoretical properties.

In the present paper, we study optimal estimation of
the mean function in two different settings. One is when the
observations are sampled at the same locations across
curves, that is, $T_{1j}=T_{2j}=\cdots=T_{nj}=:T_j$ for all
$j=1,\ldots,m$. We shall refer to this setting as \textit{common design}
because the sampling locations are common to all curves.
Another setting is when the $T_{ij}$ are independently sampled from
$\Tcal$, which we shall refer to as \textit{independent design}.
We establish the optimal rates of convergence for estimating the
mean function in both settings. Our analysis reveals interesting and different
phase transition phenomena in the two cases.
Another interesting contrast between the two settings is that
smoothing is necessary under the independent design, while, somewhat
surprisingly, it is not essential under the common design.
We remark that under the independent design, the number of sampling
points oftentimes varies from curve to curve and may even be random
itself. However, for ease of presentation and better illustration of
similarities and differences between the two types of designs, we shall
assume an equal number of sampling points on each curve in the
discussions given in this section. 

Earlier studies of nonparametric estimation of the mean function $g_0$
from a collection of discretely sampled curves can be traced back to
at least \citet{hw86} and \citet{rs91} in the
case of common design. In this setting, ignoring the temporal nature of
$\{T_j\dvtx1\le j\le m\}$, the problem of estimating $g_0$ can be
translated into estimating the mean vector $(g_0(T_1),\ldots
,g_0(T_m))'$, a typical problem in multivariate analysis. Such notions
are often quickly discarded because they essentially lead to estimating
$g_0(T_j)$ by its sample mean
%
%
\begin{equation}
\bar{Y}_{\cdot j}={1\over n}\sum_{i=1}^n Y_{ij},
\end{equation}
based on the standard Gauss--Markov theory [see, e.g., \citet{rs91}].

Note that $\bE(Y_{ij}|T)=g_0(T_{ij})$ and that the smoothness of $X$ implies
that $g_0$ is also smooth. It is therefore plausible to assume that
smoothing is essential for optimal estimation of $g_0$.
For example, a natural approach for estimating $g_0$ is to regress
$Y_{ij}$ on $T_{ij}$ nonparametrically via kernel or spline
smoothing. Various methods have been introduced along this vein [see,
e.g., \citet{rs91}]. However, not much is known about their
theoretical properties. It is noteworthy that this setting differs from
the usual nonparametric smoothing in that the observations from the
same curve are highly correlated. Nonparametric smoothing with certain
correlated errors has been previously studied by \citet{hh90},
\citet{w96} and \citet{js97}, among others.
Interested readers are referred to \citet{owy01} for a
recent survey of existing results. But neither of these earlier
developments can be applied to account for the dependency induced by
the functional nature in our setting. To comprehend the effectiveness
of smoothing in the current context, we establish minimax bounds on the
convergence rate of the integrated squared error for estimating $g_0$.

Under the common design, it is shown that the minimax rate is of the order
$m^{-2r}+n^{-1}$ where the two terms can be attributed to
discretization and stochastic error, respectively. This rate is
fundamentally different from the usual
nonparametric rate of $(nm)^{-2r/(2r+1)}$ when observations are
obtained at $nm$ distinct locations in order to recover an $r$ times
differentiable function [see, e.g., \citet{s82}]. The rate obtained
here is jointly determined by the sampling frequency
$m$ and the number of curves $n$ rather than the total number of
observations $mn$. A distinct feature of the rate is the phase
transition which occurs when $m$ is of the order $n^{1/2r}$. When the
functions are sparsely sampled, that is, $m=O(n^{1/2r})$, the optimal
rate is of the order $m^{-2r}$, solely determined by the sampling
frequency. On the other hand, when the sampling frequency is high, that is,
$m\gg n^{1/2r}$, the optimal rate remains $1/n$ regardless of
$m$. Moreover, our development uncovers a\vspace*{1pt} surprising fact that
\textit{interpolation} of $\{(T_j, \bar{Y}_{\cdot j})\dvtx j=1,\ldots,
m\}$,
that is, estimating $g_0(T_j)$ by $\bar{Y}_{\cdot j}$, is rate
optimal. In other words, contrary to the conventional wisdom, smoothing
does not result in improved convergence rates.

In addition to the common design, another popular sampling scheme is
the independent design where the $T_{ij}$ are independently sampled
from $\Tcal$.
A natural approach is to smooth observations from each curve separately
and then average over all smoothed estimates. However, the success of
this two-step procedure hinges upon the availability of a reasonable
estimate for each individual curve. In contrast to the case of common
design, we show that under the independent design, the minimax rate for
estimating $g_0$ is $(nm)^{-2r/(2r+1)}+n^{-1}$, which can be attained
by smoothing $\{(T_{ij},Y_{ij})\dvtx1\le i\le n, 1\le j\le m\}$
altogether. This implies that in the extreme case of $m=1$, the optimal
rate of estimating $g_0$ is $n^{-2r/(2r+1)}$, which also suggests the
sub-optimality of the aforementioned two-step procedure because it is
impossible to smooth a curve with only a single observation. Similar to
the common design, there is a phase transition phenomenon in the
optimal rate of convergence with a boundary at $m = n^{1/2r}$. When the
sampling frequency $m$ is small, that is, $m=O(n^{1/2r})$, the optimal
rate is of the order $(nm)^{-2r/(2r+1)}$ which depends jointly on the
values of both $m$ and $n$. In the case of high sampling frequency with
$m\gg n^{1/2r}$, the optimal rate is always $1/n$ and does not depend
on~$m$.

It is interesting to compare the minimax rates of convergence in the
two settings. The phase transition boundary for both designs occurs at
the same value, $m = n^{1/2r}$.
When $m$ is above the boundary, that is, $m \ge n^{1/2r}$, there
is no difference between the common and independent designs, and both
have the optimal rate of~$n^{-1}$. When $m$ is below the boundary,
that is, $m \ll n^{1/2r}$, the independent
design is always superior to the common design in that it offers a
faster rate of convergence. 



Our results connect with several observations made earlier in the
literature on longitudinal and functional data analysis. Many
longitudinal studies follow the independent design, and the number of
sampling points on each curve is typically small. In such settings, it
is widely recognized that one needs to pool the data to obtain good
estimates, and the two-step procedure of averaging the smooth curves
may be suboptimal. Our analysis here provides a rigorous justification
for such empirical observations by pinpointing to what extent the
two-step procedure is suboptimal. The phase transition observed here
also relates to the earlier work by \citet{hmw06} on
estimating eigenfunctions of the covariance kernel when the number of
sampling points is either fixed or of larger than $n^{1/4+\delta}$ for
some $\delta>0$. It was shown that the eigenfunctions can be estimated
at the rate of $n^{-4/5}$ in the former case and $1/n$ in the latter.
We show here that estimating the mean function has similar behavior.
Furthermore, we characterize the exact nature of such transition
behavior as the sampling frequency changes.

The rest of the paper is organized as follows. In Section
\ref{CommonDesignsec} the optimal rate of convergence
under the common design is established. We first derive a minimax lower
bound and then
show that the lower bound is in fact rate sharp. This is accomplished
by constructing a rate-optimal smoothing splines estimator. The
minimax upper bound is obtained separately for the common fixed design
and common random design. Section \ref{IndDesignsec} considers the
independent design and establishes the optimal rate of convergence in
this case. The rate-optimal estimators are easily implementable.
Numerical studies are carried out in Section \ref{Numericsec}
to demonstrate the theoretical results.
Section \ref{discussionsec} discusses connections and differences of
our results with other related work. All proofs are relegated to
Section \ref{proofsec}.

\section{Optimal rate of convergence under common design}
\label{CommonDesignsec}

In this section we consider the common design where each curve is
observed at the same set of locations $\{T_{j}\dvtx1\le j\le m\}$. We first
derive a minimax lower bound and then show that this lower bound is
sharp by constructing a smoothing splines estimator that attains the
same rate of convergence as the lower bound.

\subsection{Minimax lower bound}

Let $\Pcal(r;M_0)$ be the collection of probability measures for a
random function $X$ such that its sample path is $r$ times
differentiable almost surely and
%
%
\begin{equation}
\bE\int_\Tcal\bigl[X^{(r)}(t)\bigr]^2\,dt\le M_0
\end{equation}
for some constant $M_0>0$. Our first main result establishes the
minimax lower bound for estimating the mean function over
$\Pcal(r;M_0)$ under the common design.
\begin{theorem}
\label{thestcommonlo}
Suppose the sampling locations are common in model
(\ref{eqobsmod}). Then there exists a constant $d>0$ depending only
on $M_0$ and the variance $\sigma_0^2$ of measurement error $\varepsilon
_{ij}$ such that for any estimate $\tilde{g}$ based on observations $\{
(T_{j}, Y_{ij})\dvtx1\le i\le n, 1\le j\le m\}$,
%
%
\begin{equation}
\limsup_{n\to\infty} \sup_{\Lcal(X)\in\Pcal(r;M_0)} P\bigl(\|
\tilde{g}-g_0\|_{\Lcal_2}^2>d(m^{-2r}+n^{-1})\bigr)>0.
\end{equation}
\end{theorem}

The lower bound established in Theorem \ref{thestcommonlo} holds true
for both common fixed design where $T_j$'s are deterministic, and common
random design where $T_j$'s are also random. The term
$m^{-2r}$ in the lower bound is due to the deterministic approximation
error, and the
term $n^{-1}$ is attributed to the stochastic error. It is clear that
neither can be further improved. To see this, first consider the
situation where there is no stochastic variation and the mean
function $g_0$ is observed exactly at the points $T_j$, $j=1,\ldots,
m$. It is well known [see, e.g., \citet{dl93}] that due to
discretization, it is not possible to recover $g_0$ at a rate faster
than $m^{-2r}$ for all $g_0$ such that $\int[g_0^{(r)}]^2\le M_0$. On
the other hand, the second term $n^{-1}$ is inevitable since the mean
function $g_0$
cannot be estimated at a faster rate even if the whole random functions
$X_1,\ldots, X_n$ are observed completely. We shall show later in this
section that the rate given in the lower bound is optimal in that it
is attainable by a smoothing splines estimator.

It is interesting to notice the phase transition phenomenon in the minimax
bound. When the sampling frequency $m$ is large, it has no effect on
the rate of convergence, and $g_0$ can be estimated at the rate of
$1/n$, the best possible rate when the whole functions were
observed. More surprisingly,\vspace*{1pt} such saturation occurs when $m$ is rather
small, that is, of the order $n^{1/2r}$. On the other hand, when the
functions are sparsely sampled, that is, $m=O(n^{1/2r})$, the rate is
determined only by the sampling frequency $m$. Moreover, the rate
$m^{-2r}$ is in fact also the optimal interpolation rate. In other
words, when the functions are sparsely sampled, the mean function
$g_0$ can be estimated as well as if it is observed directly without
noise.

The rate is to be contrasted with the usual nonparametric regression
with $nm$ observations at arbitrary locations. In such a setting, it is
well known [see, e.g., \citet{t09}] that the optimal rate for
estimating $g_0$ is $(mn)^{-{2r/(2r+1)}}$, and typically stochastic
error and approximation error are of the same order to balance the
bias-variance trade-off.

\subsection{Minimax upper bound: Smoothing splines estimate}

We now consider the upper bound for the minimax risk and construct
specific rate optimal estimators under the common design. These upper
bounds show that the rate of convergence given in the lower bound
established in Theorem \ref{thestcommonlo} is sharp.
More specifically, it is shown that a smoothing splines estimator
attains the optimal rate of convergence over the parameter space
$\Pcal(r;M_0)$.

We shall\vspace*{1pt} consider a smoothing splines type of estimate suggested by
\citet{rs91}. Observe that $f \mapsto\int[f^{(r)}]^2$
is a squared semi-norm and therefore convex. By Jensen's inequality,
%
%
\begin{equation}
\int_\Tcal\bigl[g_0^{(r)}(t)\bigr]^2\,dt\le\bE\int_\Tcal
\bigl[X^{(r)}(t)\bigr]^2\,dt<\infty,
\end{equation}
which implies that $g_0$ belongs to the $r$th order Sobolev--Hilbert space,
\begin{eqnarray*}
\Wcal_2^r([0,1])&=&\bigl\{g\dvtx[0,1]\to\bR| g,g^{(1)},\ldots
,g^{(r-1)} \\
&&\hspace*{2pt}\mbox{ are absolutely
continuous and } g^{(r)}\in\Lcal_2([0,1])\bigr\}.
\end{eqnarray*}
Taking this into account, the following smoothing splines estimate can
be employed to estimate $g_0$:
%
%
\begin{equation}
\label{eqdefss0}
\hat{g}_\lambda=\argmin_{g\in\Wcal_2^r}\Biggl\{{1\over nm}\sum
_{i=1}^n\sum_{j=1}^m\bigl(Y_{ij}-g(T_{j})\bigr)^2 +\lambda\int
_\Tcal\bigl[g^{(r)}(t)\bigr]^2\,dt\Biggr\},
\end{equation}
where $\lambda>0$ is a tuning parameter that balances the fidelity to
the data and the smoothness of the estimate.

Similarly to the smoothing splines for the usual nonparametric
regression, $\hat{g}_\lambda$ can be conveniently computed, although
the minimization is taken over an infinitely-dimensional functional
space. First observe that $\hat{g}_\lambda$ can be equivalently
rewritten as
%
%
\begin{equation}
\label{eqdefss}
\hat{g}_\lambda=\argmin_{g\in\Wcal_2^r}\Biggl\{{1\over m}\sum
_{j=1}^m\bigl(\bar{Y}_{\cdot j}-g(T_{j})\bigr)^2 +\lambda\int
_\Tcal\bigl[g^{(r)}(t)\bigr]^2\,dt\Biggr\}.
\end{equation}
Appealing to the so-called representer theorem [see, e.g., \citet
{w90}], the solution of the minimization problem can be expressed as
%
%
\begin{equation}
\label{eqrep}
\hat{g}_\lambda(t)=\sum_{k=0}^{r-1}d_k t^k + \sum_{j=1}^m c_i K(t,T_j)
\end{equation}
for some coefficients $d_0,\ldots, d_{r-1},c_1,\ldots, c_m$, where
%
%
\begin{equation}
K(s,t)={1\over(r!)^2}B_r(s)B_r(t)-{1\over(2r)!}B_{2r}(|s-t|),
\end{equation}
where $B_m(\cdot)$ is the $m$th Bernoulli polynomial. Plugging
(\ref{eqrep}) back into (\ref{eqdefss}), the coefficients and
subsequently $\hat{g}_\lambda$ can be solved in a straightforward
way. This observation makes the smoothing splines procedure easily
implementable.
The readers are referred to \citet{w90} for further details.

Despite the similarity between $\hat{g}_\lambda$ and the smoothing
splines estimate in the usual nonparametric regression, they have very
different asymptotic properties. It is shown in the following that
$\hat{g}_\lambda$ achieves the lower bound established in
Theorem~\ref{thestcommonlo}. 


The analyses for the common fixed design and the common random design are
similar, and we shall focus on the fixed design where the common
sampling locations $T_1,\ldots,T_m$ are deterministic. In this case,
we assume without loss of generality that $T_1\le T_2\le\cdots\le T_m$.
The following theorem
shows that the lower bound established in Theorem \ref{thestcommonlo}
is attained by the smoothing splines estimate~$\hat{g}_\lambda$.
\begin{theorem}
\label{thsi}
Consider the common fixed design and assume that
%
%
\begin{equation}
\label{eqlocdis}
{\max_{0\le j\le m}} |T_{j+1}-T_j| \le C_0m^{-1}
\end{equation}
for some constant $C_0>0$ where we follow the convention that $T_0=0$
and \mbox{$T_{m+1}=1$}. Then
%
%
\begin{equation}
\label{eqsiconv}
\lim_{D\to\infty} \limsup_{n\to\infty} \sup_{\Lcal(X)\in\Pcal
(r;M_0)} P\bigl(\|\hat{g}_\lambda-g_0\|_{\Lcal_2}^2>D
(m^{-2r}+n^{-1})\bigr)=0
\end{equation}
for any $\lambda=O(m^{-2r}+n^{-1})$.
\end{theorem}

Together with Theorem \ref{thestcommonlo}, Theorem \ref{thsi} shows
that $\hat{g}_\lambda$ is minimax rate optimal if the tuning parameter
$\lambda$ is set to be of the order $O(m^{-2r}+n^{-1})$. We note the
necessity of the condition given by (\ref{eqlocdis}). It is clearly
satisfied when the design is equidistant, that is, $T_j=2j/(2m+1)$. The
condition ensures that the random functions are observed on a
sufficiently regular grid.

It is of conceptual importance to compare the rate of
$\hat{g}_\lambda$ with those generally achieved in the usual
nonparametric regression setting. Defined by (\ref{eqdefss0}),
$\hat{g}_\lambda$ essentially regresses $Y_{ij}$ on $T_j$. Similarly to
the usual nonparametric regression, the validity of the estimate is
driven by $\bE(Y_{ij}|T_j)=g_0(T_j)$. The difference, however, is that
$Y_{i1},\ldots, Y_{im}$ are highly correlated because they are
observed from the same random function $X_i(\cdot)$.
When all the $Y_{ij}$'s are independently\vspace*{1pt} sampled at $T_j$'s, it can be
derived that the optimal
rate for estimating $g_0$ is $m^{-2r}+(mn)^{-2r/(2r+1)}$. As we show
here, the dependency induced by the functional nature of our problem
leads to the different rate $m^{-2r}+n^{-1}$.


A distinct feature of the behavior of $\hat{g}_\lambda$ is in the
choice of the tuning parameter~$\lambda$. Tuning parameter selection
plays a paramount role in the usual nonparametric regression, as it
balances the the tradeoff between bias and variance. Optimal choice of
$\lambda$ is of the order $(mn)^{-2r/(2r+1)}$ in the usual
nonparametric regression. In contrast, in our setting, more flexibility
is allowed in the choice of the tuning parameter in that
$\hat{g}_\lambda$ is rate optimal so long as $\lambda$ is sufficiently
small. In particular, taking $\lambda\to0^+$, $\hat{g}_\lambda$ reduces
to the splines interpolation, that is, the solution to
%
%
\begin{equation}
\label{eqinterpolation}\quad
\min_{g\in\Wcal_2^r} \int_\Tcal\bigl[g^{(r)}(t)\bigr]^2
\qquad\mbox{subject to } g(T_j)=\bar{Y}_{\cdot j},\qquad j=1,\ldots, m.
\end{equation}
This amounts to, in particular, estimating $g_0(T_j)$ by
$\bar{Y}_{\cdot j}$. In other words, there is no benefit from
smoothing in terms of the convergence rate. However, as we will see in
Section \ref{Numericsec}, smoothing can lead to improved finite
sample performance.
\begin{Remark*} More general statements can also be made without the
condition on the spacing of sampling points. More specifically, denote
by
\[
R(T_1,\ldots,T_m)=\max_j |T_{j+1}-T_j|
\]
the discretization resolution. Using the same argument, one can show
that the optimal convergence rate in the minimax sense is
$R^{2r}+n^{-1}$ and $\hat{g}_\lambda$ is rate optimal so long as
$\lambda=O(R^{2r}+n^{-1})$.
\end{Remark*}
\begin{Remark*}
Although we have focused here on the case when the
sampling points are deterministic, a similar statement can also be made
for the setting where the sampling points are random. In particular,
assuming that $T_j$ are independent and identically distributed with a
density function $\eta$ such that $\inf_{t\in\Tcal} \eta(t)\ge
c_0>0$ and $g_0\in\Wcal_\infty^r$, it can be shown that the
smoothing splines estimator $\hat g_\lambda$ satisfies
%
%
\begin{equation}\quad
\label{eqranup}
\lim_{D\to\infty} \limsup_{n\to\infty} \sup_{\Lcal(X)\in\Pcal
(r; M_0)}
P\bigl(\|\hat{g}_\lambda-g_0\|_{\Lcal_2}^2>D
(m^{-2r}+n^{-1})\bigr)=0
\end{equation}
for any $\lambda=O(m^{-2r}+n^{-1})$. In other words, $\hat{g}_\lambda
$ remains rate optimal. 
\end{Remark*}




\section{Optimal rate of convergence under independent design}
\label{IndDesignsec}

In many applications, the random functions $X_i$ are not observed at common
locations. Instead, each curve is discretely observed at a different
set of points [see, e.g., \citet{jh01}, \citet{rw01},
\citet{dhlz02}, \citet{ymw05}]. In these
settings, it is more
appropriate to model the sampling points $T_{ij}$ as independently
sampled from a common distribution. In this section we shall consider
optimal estimation of the mean function under the
independent design.

Interestingly, the behavior of the estimation problem is drastically
different between the common design and the independent design. To keep
our treatment general, we allow the number of sampling points to vary.
Let
$m$ be the harmonic mean of $m_1,\ldots,m_n$, that is,
\[
m:=\Biggl({1\over n}\sum_{i=1}^n {1\over m_i}\Biggr)^{-1}.
\]
Denote by $\Mcal(m)$ the collection of sampling frequencies
$(m_1,\ldots,m_n)$ whose harmonic mean is $m$.
In
parallel to Theorem \ref{thestcommonlo}, we have the following
minimax lower bound for estimating $g_0$ under the independent design.
\begin{theorem}
\label{thestmeanlo}
Suppose $T_{ij}$ are independent and identically distributed with a
density function $\eta$ such that $\inf_{t\in\Tcal} \eta(t)\ge
c_0>0$. Then there\vspace*{1pt} exists a constant $d>0$ depending only on $M_0$ and
$\sigma_0^2$ such that for any estimate $\tilde{g}$ based on
observations $\{(T_{ij}, Y_{ij})\dvtx1\le i\le n, 1\le j\le m\}$,
%
%
\begin{equation}
\label{eqgconvlo}\qquad
\limsup_{n\to\infty}
\mathop{\sup_{\Lcal(X)\in\Pcal(r;M_0)}}_{(m_1,\ldots,m_n)\in
\Mcal(m)} P\bigl(\|\tilde{g}-g_0\|_{\Lcal_2}^2>d\bigl(
(nm)^{-{2r/(2r+1)}}+n^{-1}\bigr)\bigr)>0.
\end{equation}
\end{theorem}

The minimax lower bound given in Theorem \ref{thestmeanlo} can also
be achieved using the smoothing splines type of estimate. To account
for the different sampling frequency for different curves, we consider
the following estimate of $g_0$:
%
%
\begin{equation}
\hat{g}_\lambda=\argmin_{g\in\Wcal_2^r}\Biggl\{{1\over n}\sum
_{i=1}^n{1\over m_i}\sum_{j=1}^{m_i}\bigl(Y_{ij}-g(T_{ij})\bigr)^2
+\lambda\int_\Tcal\bigl[g^{(r)}(t)\bigr]^2\,dt\Biggr\}.
\end{equation}

\begin{theorem}
\label{thestmean}
Under the conditions of Theorem \ref{thestmeanlo},
if $\lambda\asymp(nm)^{-2r/(2r+1)}$, then the smoothing splines
estimator $\hat g_\lambda$ satisfies
%
%
\begin{eqnarray}
\label{eqgconv}\quad
&&\lim_{D\to\infty} \limsup_{n\to\infty}
\mathop{\sup_{\Lcal(X)\in\Pcal(r;M_0)}}_{(m_1,\ldots,m_n)\in
\Mcal(m)} P\bigl(\|\hat{g}-g_0\|_{\Lcal_2}^2>D\bigl(
(nm)^{-{2r/(2r+1)}}+n^{-1}\bigr)\bigr)\nonumber\\[-8pt]\\[-8pt]
&&\qquad=0.\nonumber
\end{eqnarray}
In other words, $\hat{g}_\lambda$ is rate optimal.
\end{theorem}

Theorems \ref{thestmeanlo} and \ref{thestmean} demonstrate both
similarities and significant differences between the two types of
designs in terms of the convergence rate. For either the common design
or the independent design, the sampling frequency only plays a role in
determining\vspace*{1pt} the convergence rate when the functions are sparsely
sampled, that is, $m=O(n^{1/2r})$. But how the sampling frequency affects
the convergence rate when each curve is sparsely sampled differs
between the two designs. For the independent design, the total number
of observations $mn$, whereas for the common design $m$ alone,
determines\vspace*{1pt} the minimax rate. It is also noteworthy that when
$m=O(n^{1/2r})$, the optimal rate under the independent design,
$(mn)^{-2r/(2r+1)}$, is the same as if all the observations are
independently observed. In other words, the dependency among
$Y_{i1},\ldots, Y_{im}$ does not affect the convergence rate in this case.
\begin{Remark*}
We emphasize that Theorems \ref{thestmeanlo} and \ref
{thestmean} apply to both deterministic and random sampling
frequencies. In particular for random sampling frequencies, together
with the law of large numbers, the same minimax bound holds when we
replace the harmonic mean by
\[
\Biggl({1\over n}\sum_{i=1}^n\bE(1/m_i)\Biggr)^{-1},
\]
when assuming that $m_i$'s are independent.
\end{Remark*}

\subsection{Comparison with two-stage estimate}

A popular strategy to handle discretely sampled functional data in
practice is a two-stage procedure. In the first step, nonparametric
regression is run for data from each curve to obtain estimate $\tilde
{X}_i$ of $X_i$, $i=1,2,\ldots,n$. For example, they can be obtained
by smoothing splines
%
%
\begin{equation}
\tilde{X}_{i,\lambda}=\argmin_{f\in\Wcal_2^r}\Biggl\{{1\over
m}\sum_{j=1}^m\bigl(Y_{ij}-g(T_{ij})\bigr)^2 +\lambda\int_\Tcal
\bigl[f^{(r)}(t)\bigr]^2\,dt\Biggr\}.
\end{equation}
Any subsequent inference can be carried out using the
$\tilde{X}_{i,\lambda}$ as if they were the original true random
functions. In particular, the mean function $g_0$ can be estimated by
the simple average
%
%
\begin{equation}
\tilde{g}_\lambda={1\over n}\sum_{i=1}^n \tilde{X}_{i,\lambda}.
\end{equation}
Although formulated differently, it is worth pointing out that this
procedure is equivalent to the smoothing splines estimate $\hat
{g}_\lambda$ under the common design.
\begin{proposition}
\label{prequiv}
Under the common design, that is, $T_{ij}=T_j$ for $1\le i\le n$ and
$1\le j\le m$. The estimate $\tilde{g}_\lambda$ from the two-stage
procedure is equivalent to the smoothing splines estimate $\hat
{g}_\lambda$: $\hat{g}_\lambda=\tilde{g}_\lambda$.
\end{proposition}

In light of Theorems \ref{thestcommonlo} and \ref{thsi}
, the two-step procedure is also rate optimal under the common design.
But it
is of great practical importance to note that in order to achieve the
optimality, it is critical that in the first step we \textit{undersmooth}
each curve by using a sufficiently small tuning parameter.

Under independent design, however, the equivalence no long holds. The
success of the two-step estimate $\tilde{g}_\lambda$ depends upon
getting a good estimate of each curve, which is not possible when $m$
is very small. In the extreme case of $m=1$, the procedure is no longer
applicable, but Theorem \ref{thestmean} indicates\vspace*{1pt} that
smoothing splines estimate $\hat{g}_\lambda$ can still achieve the
optimal convergence rate of $n^{-2r/(2r+1)}$. The readers are also
referred to \citet{hmw06} for discussions on the pros and cons of
similar two-step procedures in the context of estimating the functional
principal components.

\section{Numerical experiments}
\label{Numericsec}

The smoothing splines estimators are easy to implement.
To demonstrate the practical implications of our theoretical results,
we carried out a set of simulation studies. The true mean function
$g_0$ is fixed as
%
%
\begin{equation}
g_0=\sum_{k=1}^{50} 4(-1)^{k+1}k^{-2} \phi_k,
\end{equation}
where $\phi_1(t)=1$ and $\phi_{k+1}(t)=\sqrt{2}\cos(k\pi t)$ for
$k\ge1$. The random function $X$ was generated as
%
%
\begin{equation}
X=g_0+\sum_{k=1}^{50} \zeta_kZ_k \phi_k,
\end{equation}
where $Z_k$ are independently sampled from the uniform distribution
on\break
$[-\sqrt{3},\sqrt{3}]$, and $\zeta_k$ are deterministic. It is not
hard to see that $\zeta_k^2$ are the eigenvalues of the covariance
function of $X$ and therefore determine the smoothness of a sample
curve. In particular, we take $\zeta_k=(-1)^{k+1}k^{-1.1/2}$. It is
clear that the sample path of $X$ belongs to the second order Sobolev
space ($r=2$).

%
\begin{figure}[b]

\includegraphics{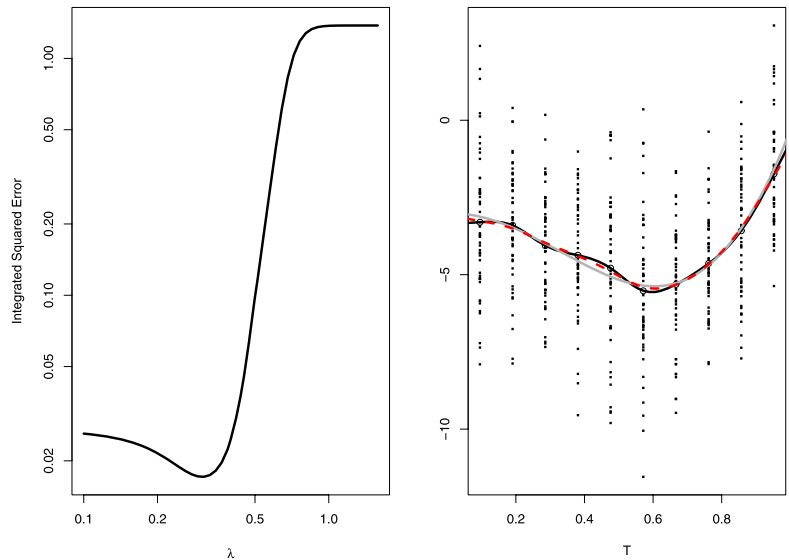}

\caption{Effect of smoothing under common design: for a typical data
set with fifty curves, ten observations were taken on each curve. The
observations and $g_0$ (solid grey line) are given in the right panel
together with the spline interpolation estimate (solid black line) and
smoothing splines estimate (red dashed line) with the tuning parameter
chosen to yield the smallest integrated squared error. The left panel
gives the integrated squared error of the smoothing splines estimate as
a function of the tuning parameter. It is noted that the smoothing
splines estimate essentially reduces to the spline interpolation for
$\lambda$ smaller than $0.1$.}
\label{figssi}
\end{figure}

We begin with a set of simulations designed to demonstrate the effect
of interpolation and smoothing under common design. A data set of fifty
curves were first simulated according to the aforementioned scheme. For
each curve, ten noisy observations were taken at equidistant locations
on each curve following model (\ref{eqobsmod}) with $\sigma
_0^2=0.5^2$. The observations, together with $g_0$ (grey line), are
given in the right panel of Figure \ref{figssi}. Smoothing splines
estimate $\hat{g}_\lambda$ is also computed with a variety of values
for $\lambda$. The integrated squared error, $\|\hat{g}_\lambda-g_0\|
_{\Lcal_2}$, as a function of the tuning parameter $\lambda$ is given
in the left panel. For $\lambda$ smaller than $0.1$, the smoothing
splines estimate essentially reduces to the spline interpolation. To
contrast the effect of interpolation and smoothing, the right panel
also includes the interpolation estimate (solid black line) and $\hat
{g}_\lambda$ (red dashed line) with the tuning parameter chosen to
minimize the integrated squared error. We observe from the figure that
smoothing does lead to slightly improved finite sample performance
although it does not affect the convergence rate as shown in Section
\ref{CommonDesignsec}.

\begin{figure}[b]

\includegraphics{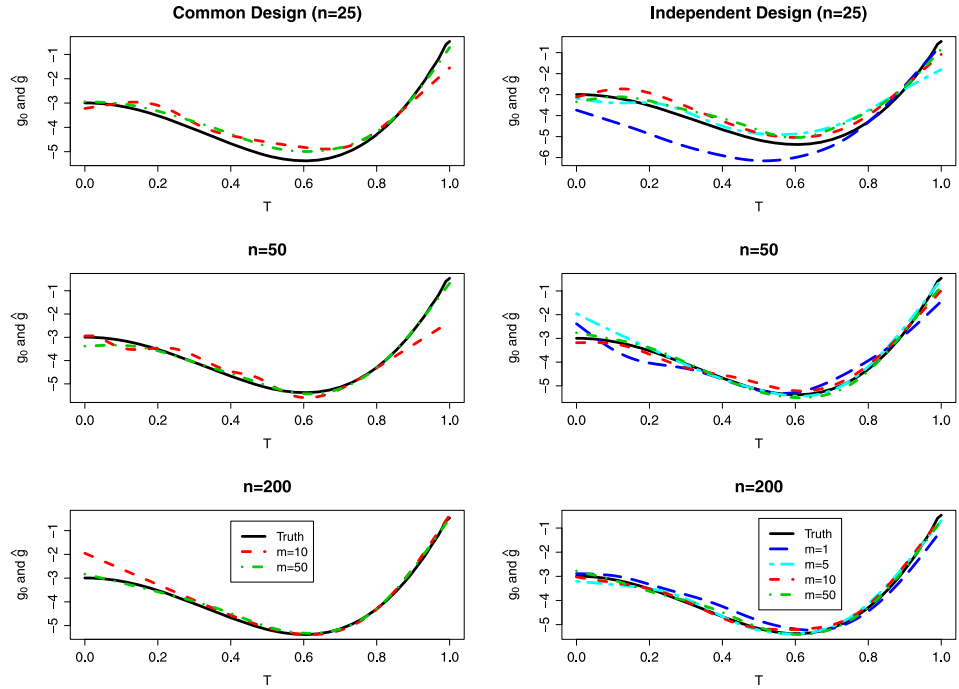}

\caption{Effect of $m$, $n$ and type of design on estimating $g_0$:
smoothing splines estimates obtained under various combinations are
plotted together with $g_0$.} \label{figtypical}
\end{figure}

The next numerical experiment intends to demonstrate the effect of sample
size~$n$, sampling frequency $m$ as well as design. To this end, we
simulated $n$ curves, and from each curve, $m$ discrete observations
were taken following model (\ref{eqobsmod}) with
$\sigma_0^2=0.5^2$. The sampling locations are either fixed at
$T_j=(2j)/(2m+1)$, $j=1,\ldots, m$, for common design or randomly
sampled from the uniform distribution on $[0,1]$. The smoothing
splines estimate $\hat{g}_\lambda$ for each simulated data set, and the
tuning parameter is set to yield the smallest integrated squared error
and therefore reflect the best performance of the estimating procedure
for each data set. We repeat the experiment with varying combinations
of $n=25, 50$ or $200$, $m=1,5, 10$ or $50$. For the common design, we
restrict to $m=10$ or $50$ to give more meaningful comparison. The true
function $g_0$ as well as its estimates obtained in each of the
settings are given in Figure \ref{figtypical}.

%

Figure \ref{figtypical} agrees pretty well with our theoretical
results. For instance, increasing either $m$ or $n$ leads to improved
estimates, whereas such improvement is more visible for small values of
$m$. Moreover, for the same value of $m$ and $n$, independent designs
tend to yield better estimates.

\begin{figure}

\includegraphics{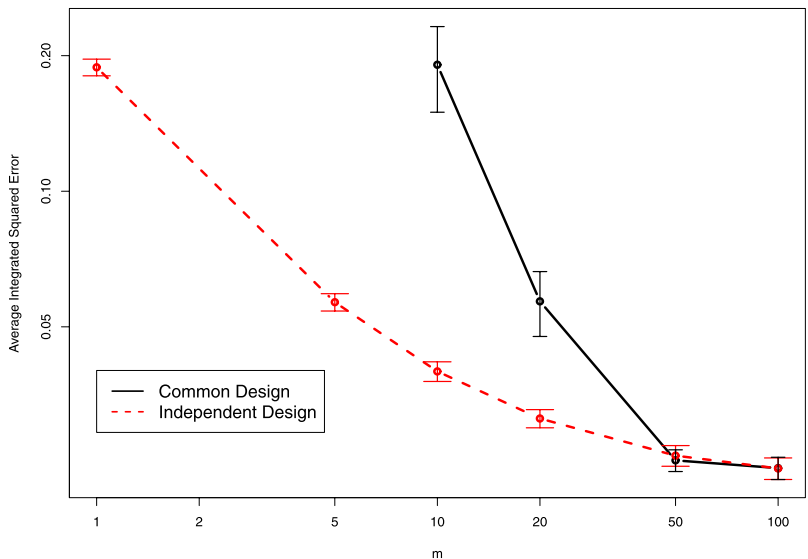}

\caption{Effect of design type and sampling frequency on estimating
$g_0$: the black solid line and circles correspond to common design
whereas the red dashed lines and circles correspond to independent
design. The error bars correspond to the average $\pm$ one standard
errors based on two hundred repetitions. Note that both axes are in log
scale to yield better comparison.}
\label{figdesignfig}
\end{figure}

To further contrast the two types of designs and the effect of
sampling frequency on estimating $g_0$, we now fix the number of
curves at $n=100$. For the common design, we consider $m=10, 20, 50$
or $100$. For the independent design, we let $m=1, 5, 10, 20, 50$ or
$100$. For each combination of $(n,m)$, two hundred data sets were
simulated following the same mechanism as before. Figure
\ref{figdesignfig} gives the estimation error averaged over the one
hundred data sets for each combination of $(n,m)$. It clearly shows
%
%
that independent design is preferable over common design when $m$ is
small; and the two types of designs are similar when $m$ is
large. Both phenomena are in agreement with our theoretical results
developed in the earlier sections.

\section{Discussions}
\label{discussionsec}
We have established the optimal rates of convergence for estimating
the mean function under both the common design and independent
design. The results reveal several significant differences in the
behavior of the minimax estimation problem between the two designs.
These revelations have important theoretical and practical implications.
In particular, for sparsely sampled functions,
the independent design leads to a faster rate of convergence
when compared to the common design and thus should be
preferred in practice.


The optimal rates of convergence for estimating the mean function based
on discretely sampled random functions behave in a fundamentally
different way from the minimax rate of convergence in the conventional
nonparametric regression problems. The optimal rates in the mean
function estimation are jointly determined by the sampling frequency
$m$ and the number of curves $n$ rather than the total number of
observations $mn$.

The observation that one can estimate the mean function as well as if
the the whole curves are available when $m\gg n^{1/2r}$ bears some
similarity to some recent findings on estimating the covariance kernel
and its eigenfunction under independent design. Assuming that $X$ is
twice differentiable (i.e., $r=2$), \citet{hmw06}
showed that when $m\gg n^{1/4+\delta}$ for some $\delta>0$, the
covariance kernel and its eigenfunctions can be estimated at the rate
of $1/n$ when using\vspace*{1pt} a two-step procedure. More recently, the cutoff
point is further improved to $n^{1/2r}\log n$ with general $r$ by
\citet{cy10b} using an alternative method. Intuitively one may expect
estimating the covariance kernel to be more difficult than estimating
the mean function, which suggests that these results may not be
improved much further for estimating the covariance kernel or its
eigenfunctions.

We have also shown that the particular smoothing splines type of
estimate discussed earlier by \citet{rs91} attains the
optimal convergence rates under both designs with appropriate
tuning. The smoothing splines estimator is well suited for nonparametric
estimation over Sobolev spaces.
We note, however, other nonparametric techniques such as
kernel or local polynomial estimators can also be used. We expect that
kernel smoothing or other methods with proper choice of the tuning
parameters can also achieve the optimal rate of convergence. Further
study in this direction is beyond the scope of the current paper, and
we leave it for future research.


Finally, we emphasize that although we have focused on the univariate
Sobolev space for simplicity, the phenomena observed and techniques
developed apply to more general functional spaces. Consider, for
example, the multivariate setting where $\Tcal=[0,1]^d$. Following
the same arguments, it can be shown that the minimax rate for
estimating an $r$-times differentiable function is $m^{-2r/d}+n^{-1}$
under the common design and $(nm)^{-2r/(2r+d)}+n^{-1}$ under the
independent design. The phase transition phenomena thus remain in the
multidimensional setting under both designs with a transition boundary
of $m=n^{d/2r}$.

\section{Proofs}
\label{proofsec}

\mbox{}

\begin{pf*}{Proof of Theorem \ref{thestcommonlo}}
Let $\Dcal$ be the collection all measurable functions of $\{
(T_{ij},Y_{ij})\dvtx1\le i\le n, 1\le j\le m\}$. First note that it is
straightforward to show that
\[
\limsup_{n\to\infty} \inf_{\tilde{g}\in\Dcal}\sup_{\Lcal(X)\in
\Pcal(r;M_0)} P(\|\tilde{g}-g_0\|_{\Lcal_2}^2>dn^{-1})>0
\]
by considering $X$ as an unknown constant function where the problem
essentially becomes estimating the mean from $n$ i.i.d. observations,
and $1/n$ is known as the optimal rate. It now suffices to show that
\[
\limsup_{n\to\infty} \inf_{\tilde{g}\in\Dcal}\sup_{\Lcal(X)\in
\Pcal(r;M_0)} P(\|\tilde{g}-g_0\|_{\Lcal_2}^2>dm^{-2r})>0.
\]
Let $\varphi_1,\ldots,\varphi_{2m}$ be $2m$ functions from $\Wcal
_2^r$ with distinct support, that is,
\[
\varphi_k(\cdot)=h^rK\biggl({\cdot-t_k\over h}\biggr),\qquad k =1,\ldots, 2m,
\]
where $h=1/(2m)$, $t_k=(k-1)/2m+1/4m$, and $K\dvtx\bR\to[0,\infty)$ is
an $r$ times differentiable function with support $[-1/2, 1/2]$. See
\citet{t09} for explicit construction of such functions.

For each $b=(b_1,\ldots, b_{2m})\in\{0,1\}^{2m}$, define
\[
g_b(\cdot)=\sum_{k=1}^{2m} b_{k}\varphi_k(\cdot).
\]
It is clear that
\[
\min_{H(b,b')\ge1}{\|g_b-g_{b'}\|^2_{\Lcal_2}\over H(b,b')}=\|
\varphi_k\|_{\Lcal_2}^2=(2m)^{-(2r+1)}\|K\|_{\Lcal_2}^2.
\]
The claim then follows from an application of Assouad's lemma
[\citet{a83}].
\end{pf*}
\begin{pf*}{Proof of Theorem \ref{thsi}}
It is well known [see, e.g., \citet{gs94}] that $\hat
{g}_\lambda$ can be characterized as the solution to the following:
\[
\min_{g\in\Wcal_2^r} \int_\Tcal\bigl[g^{(r)}(t)\bigr]^2\,dt
\qquad\mbox{subject to } g(T_j)=\hat{g}_\lambda(T_j),\qquad
j=1,\ldots, m.
\]
Write
\[
{\delta}_{j}=\hat{g}_\lambda(T_{ij})-g_0(T_{ij}),
\]
and let $h$ be the linear interpolation of $\{(T_j,{\delta}_{j})\dvtx
1\le
j\le m\}$, that is,
%
%
\begin{equation}
\label{eqdefh}
h(t)=\cases{
{\delta}_{1}, &\quad$0\le t\le T_1$,\vspace*{2pt}\cr
\displaystyle {\delta}_{j} {T_{j+1}-t\over T_{j+1}-T_j}+{\delta}_{\cdot j+1}
{t-T_j\over T_{j+1}-T_j}, &\quad$T_j\le t\le T_{j+1}$,\vspace*{2pt}\cr
{\delta}_{m}, &\quad $T_m\le t\le1$.}
\end{equation}
Then $\hat{g}_\lambda=Q_T(g_0+h)$ where $Q_T$ be the operator
associated with the $r$th order spline interpolation, that is, $Q_T(f)$
is the solution to
\[
\min_{g\in\Wcal_2^r} \int_\Tcal\bigl[g^{(r)}(t)\bigr]^2\,dt\qquad
\mbox{subject to } g(T_j)=f(T_j),\qquad j=1,\ldots, m.
\]

Recall that $Q_T$ is a linear operator in that
$Q_T(f_1+f_2)=Q_T(f_1)+Q_T(f_2)$ [see, e.g., \citet{dl93}].
Therefore, $\hat{g}_\lambda=Q_T(g_0)+Q_T(h)$. By the triangular inequality,
%
%
\begin{equation}
\label{eqinterproof1}
\|\hat{g}-g_0\|_{\Lcal_2}\le\|Q_T(g_0)-g_0\|_{\Lcal_2}+\|Q_T(h)\|
_{\Lcal_2}.
\end{equation}
The first term on the right-hand side represents the approximation
error of spline interpolation for $g_0$, and it is well known that it
can be bounded by [see, e.g., \citet{dl93}]
%
%
\begin{eqnarray}
\label{eqinterproof2}
&&\|Q_T(g_0)-g_0\|_{\Lcal_2}^2\nonumber\\
&&\qquad\le c_0\Bigl(\max_{0\le j\le
m}|T_{j+1}-T_j|^{2r} \Bigr)\int_\Tcal\bigl[g_0^{(r)}(t)
\bigr]^2\,dt\\
&&\qquad\le c_0M_0m^{-2r}.\nonumber
\end{eqnarray}
Hereafter, we shall use $c_0>0$ as a generic constant which may take
different values at different appearance.

It now remains to bound $\|Q_T(h)\|_{\Lcal_2}$. We appeal to the
relationship between spline interpolation and the best local polynomial
approximation. Let
\[
I_j=[T_{j-r+1},T_{j+r}]
\]
with the convention that $T_j=0$ for $j<1$ and $T_j=1$ for $j>m$.
Denote by $P_j$ the best approximation error that can be achieved on
$I_j$ by a polynomial of order less that $r$, that is,
%
%
\begin{equation}
\label{eqpolyappr}
P_j(f)=\min_{a_k\dvtx k<r}\int_{I_j}\Biggl[\sum
_{k=0}^{r-1}a_kt^k-f(t)\Biggr]^2\,dt.
\end{equation}
It can be shown [see, e.g., Theorem 4.5 on page 147 of \citet
{dl93}] that
\[
\|f-Q_T(f)\|_{\Lcal_2}^2\le c_0\Biggl(\sum_{j=1}^mP_j(f)^2\Biggr).
\]
Then
\[
\|Q_T(h)\|_{\Lcal_2}\le\|h\|_{\Lcal_2}+c_0\Biggl(\sum
_{j=1}^mP_j(f)^2\Biggr)^{1/2}.
\]
Together with the fact that
\[
P_j(h)^2\le\int_{I_j} h(t)^2\,dt,
\]
we have
\begin{eqnarray*}
\|Q_T(h)\|^2_{\Lcal_2}&\le& c_0\|h\|^2_{\Lcal_2}\\
&\le& c_0\sum
_{j=1}^m\delta_j^2(T_{j+1}-T_{j-1})\le c_0m^{-1}\sum_{j=1}^m\delta
_j^2\\
&=&c_0m^{-1}\sum_{j=1}^m[\hat{g}_\lambda(T_j)-g_0(T_j)
]^2\\
&\le&c_0m^{-1}\sum_{j=1}^m\bigl([\bar{Y}_{\cdot j}-\hat
{g}_\lambda(T_j)]^2+[\bar{Y}_{\cdot j}-g_0(T_j)
]^2\bigr).
\end{eqnarray*}

Observe that
\[
\bE\Biggl(m^{-1}\sum_{j=1}^m[\bar{Y}_{\cdot j}-g_0(T_j)
]^2\Biggr)=c_0\sigma_0^2n^{-1}.
\]
It suffices to show that
\[
m^{-1}\sum_{j=1}^m[\bar{Y}_{\cdot j}-\hat{g}_\lambda
(T_j)]^2=O_p(m^{-2r}+n^{-1}).
\]
To this end, note that by the definition of $\hat{g}_\lambda$,
\begin{eqnarray*}
{1\over m}\sum_{j=1}^m\bigl(\bar{Y}_{\cdot j}-\hat{g}_\lambda
(T_{j})\bigr)^2&\le& {1\over m}\sum_{j=1}^m\bigl(\bar{Y}_{\cdot
j}-\hat{g}_\lambda(T_{j})\bigr)^2 +\lambda\int_\Tcal\bigl[\hat
{g}_\lambda^{(r)}(t)\bigr]^2\,dt\\
&\le& {1\over m}\sum_{j=1}^m\bigl(\bar{Y}_{\cdot
j}-g_0(T_{j})\bigr)^2 +\lambda\int_\Tcal\bigl[g_0^{(r)}(t)
\bigr]^2\,dt\\
&\le& O_p(m^{-2r}+n^{-1}),
\end{eqnarray*}
because $\lambda=O(m^{-2r}+n^{-1})$. The proof is now complete.
\end{pf*}

\begin{pf*}{Proof of Theorem \ref{thestmeanlo}}
Note that any lower bound for a specific case yields immediately a
lower bound for the general case. It therefore suffices to consider the
case when $X$ is a Gaussian process and $m_1=m_2=\cdots=m_n=:m$.
Denote by $N=c(nm)^{1/(2r+1)}$ where $c>0$ is a constant to be
specified later. Let $b=(b_1,\ldots, b_N)\in\{0,1\}^N$ be a binary
sequence, and write
\[
g_b(\cdot)=M_0^{1/2}\pi^{-r}\sum_{k=N+1}^{2N}
N^{-1/2}k^{-r}b_{k-N}\varphi_k(\cdot),
\]
where $\varphi_k(t)=\sqrt{2}\cos(\pi kt)$. It is not hard to see that
\begin{eqnarray*}
\int_{\Tcal}\bigl[g_b^{(r)}(t)\bigr]^2\,dt&=&M_0\pi^{-2r}\sum_{k\ge
N+1}^{2N}(\pi k)^{2r}(N^{-1/2}k^{-r}b_{k-N}
)^2\\
&=&M_0N^{-1}\sum_{k= N+1}^{2N}b_{k-N}\le M_0.
\end{eqnarray*}
Furthermore,
\begin{eqnarray*}
\|g_b-g_{b'}\|^2_{\Lcal_2}&=&M_0\pi^{-2r}N^{-1}\sum
_{k=N+1}^{2N}k^{-2r}(b_{k-N}-b_{k-N}')^2\\
&\ge& M_0\pi^{-2r}(2N)^{-(2r+1)}\sum
_{k=N+1}^{2N}(b_{k-m}-b_{k-m}')^2\\
&=&c_0N^{-(2r+1)}H(b,b')
\end{eqnarray*}
for some\vspace*{1pt} constant $c_0>0$. By the Varshamov--Gilbert bound [see, e.g.,
\citet{t09}], there exists a collection of binary sequences $\{
b^{(1)},\ldots, b^{(M)}\}\subset\{0,1\}^N$ such that $M\ge2^{N/8}$, and
\[
H\bigl(b^{(j)},b^{(k)}\bigr)\ge N/8\qquad \forall1\le j<k\le M.
\]
Then
\[
\bigl\|g_{b^{(j)}}-g_{b^{(k)}}\bigr\|_{\Lcal_2}\ge c_0N^{-r}.
\]

Assume that $X$ is a Gaussian process with mean $g_b$, $T$ follows a
uniform distribution on $\Tcal$ and the measurement error $\varepsilon
\sim N(0,\sigma_0^2)$. Conditional on $\{T_{ij}\dvtx j=1,\ldots, m\}$,
$Z_{i\cdot}=(Z_{i1},\ldots, Z_{im})'$ follows a multivariate normal
distribution with mean $\mu_b=(g_b(T_{i1}),\ldots, g_b(T_{im}))'$ and
covariance\vspace*{1pt} matrix $\Sigma(T)=(C_0(T_{ij}, T_{ik}))_{1\le j,k\le
m}+\sigma^2_0 I$. Therefore, the Kullback--Leibler distance from
probability measure $\Pi_{g_{b^{(j)}}}$ to $\Pi_{g_{b^{(k)}}}$ can be
bounded by
\begin{eqnarray*}
\operatorname{KL}\bigl(\Pi_{g_{b^{(j)}}}|\Pi_{g_{b^{(k)}}}\bigr)&=&n\bE_{T} \bigl[\bigl(\mu
_{g_{b^{(j)}}}-\mu_{g_{b^{(k)}}}\bigr)'\Sigma^{-1}(T)\bigl(\mu
_{g_{b^{(j)}}}-\mu_{g_{b^{(k)}}}\bigr)\bigr]\\
&\le&n\sigma_0^{-2}\bE_{T} \bigl\|\mu_{g_{b^{(j)}}}-\mu_{g_{b^{(k)}}}\bigr\|
^2\\
&=&nm\sigma_0^{-2}\bigl\|g_{b^{(j)}}-g_{b^{(k)}}\bigr\|_{\Lcal_2}^2\\
&\le&c_1nm\sigma_0^{-2}N^{-2r}.
\end{eqnarray*}
An application of Fano's lemma now yields
\begin{eqnarray*}
\max_{1\le j\le M}\bE_{g_{b^{(j)}}}\bigl\|\tilde{g}-g_{b^{(j)}}\bigr\|_{\Lcal
_2}&\ge& c_0N^{-r}\biggl(1-{\log(c_1nm\sigma_0^{-2}N^{-2r})+\log2\over
\log M}\biggr)\\
&\asymp&(nm)^{-{r/(2r+1)}}
\end{eqnarray*}
with an appropriate choice of $c$, for any estimate $\tilde{g}$. This
in turn implies that
\[
\limsup_{n\to\infty} \inf_{\tilde{g}\in\Dcal}\sup_{\Lcal(X)\in
\Pcal(r;M_0)} P\bigl(\|\tilde{g}-g_0\|_{\Lcal_2}^2>d(nm
)^{-{2r/(2r+1)}}\bigr)>0.
\]
The proof can then be completed by considering $X$ as an unknown
constant function.
\end{pf*}
\begin{pf*}{Proof of Theorem \ref{thestmean}}
For brevity, in what follows, we treat the sampling frequencies
$m_1,\ldots, m_n$ as deterministic. All the arguments, however, also
apply to the situation when they are random by treating all the
expectations and probabilities as conditional on $m_1,\ldots,m_n$.
Similarly, we shall also assume that $T_j$'s follow uniform
distribution. The argument can be easily applied to handle more general
distributions.

It is well known that $\Wcal_2^r$, endowed with the norm
%
%
\begin{equation}
\label{eqsobnorm}
\Vert f \Vert^2_{\Wcal_2^r}=\int f^2+\int
\bigl(f^{(r)}\bigr)^2,
\end{equation}
forms a reproducing kernel Hilbert space [\citet{a50}].
Let $\Hcal_0$ be the collection of all polynomials of order less than
$r$ and $\Hcal_1$ be its orthogonal complement in $\Wcal_2^r$. Let $\{
\phi_k\dvtx1\le k\le r\}$ be a set of orthonormal basis functions of
$\Hcal_0$, and $\{\phi_k\dvtx k>r\}$ an orthonormal basis of $\Hcal_1$
such that any $f\in\Wcal_2^r$ admits the representation
\[
f=\sum_{\nu\ge1} f_{\nu}\phi_\nu.
\]
Furthermore,
\[
\|f\|_{\Lcal_2}^2=\sum_{\nu\ge1} f_{\nu}^2\quad \mbox{and}\quad
\|f\|_{\Wcal_2^r}^2=\sum_{\nu\ge1} (1+\rho_\nu^{-1})f_{\nu}^2,
\]
where $\rho_1=\cdots=\rho_r=+\infty$ and $\rho_\nu\asymp\nu^{-2r}$.


Recall that
\[
\hat{g}=\argmin_{g\in\Hcal(K)}\biggl\{\ell_{mn}(g) +\lambda\int
\bigl[g^{(r)}\bigr]^2\biggr\},
\]
where
\[
\ell_{mn}(g)={1\over n}\sum_{i=1}^n{1\over m_i}\sum_{j=1}^{m_i}
\bigl(Y_{ij}-g(T_{ij})\bigr)^2.
\]
For brevity, we shall abbreviate the subscript of $\hat{g}$ hereafter
when no confusion occurs. Write
\begin{eqnarray*}
\ell_{\infty}(g)&=&\bE\Biggl({1\over n}\sum_{i=1}^n{1\over m_i}\sum
_{j=1}^{m_i}[Y_{ij}-g(T_{ij})]^2\Biggr)\\
&=&\bE\bigl([Y_{11}-g_0(T_{11})]^2\bigr)+\int_{\Tcal
}[g(s)-g_0(s)]^2\,ds.
\end{eqnarray*}
Let
\[
\bar{g}=\argmin_{g\in\Hcal(K)}\biggl\{\ell_{\infty}(g) +\lambda
\int\bigl[g^{(r)}\bigr]^2\biggr\}.
\]
Denote
\[
\ell_{mn,\lambda}(g)=\ell_{mn}(g)+\lambda\int\bigl[g^{(r)}\bigr]^2;\qquad
\ell_{\infty,\lambda}(g)=\ell_\infty(g)+\lambda\int\bigl[g^{(r)}\bigr]^2.
\]

Let
\[
\tilde{g}=\bar{g}-\tfrac{1}{2}G_\lambda^{-1}D\ell_{mn,\lambda}(\bar{g}),
\]
where $G_\lambda= (1/2)D^2\ell_{\infty,\lambda}(\bar{g})$ and $D$
stands for the Fr\'echet derivative. It is clear that
\[
\hat{g}-g_0=(\bar{g}-g_0)+(\hat{g}-\tilde
{g})+(\tilde{g}-\bar{g}).
\]
We proceed by bounding the three terms on the right-hand side
separately. In particular, it can be shown that
%
%
\begin{equation}
\label{eqerr1}
\|\bar{g}-g_0\|^2_{\Lcal_2}\le c_0\lambda\int\bigl[g_0^{(r)}\bigr]^2
\end{equation}
and
%
%
\begin{equation}
\label{eqerr2}
\Vert\tilde{g}-\bar{g}\Vert_{\Lcal_2}^2=O_p
\bigl(n^{-1}+(nm)^{-1}\lambda^{-{1/(2r)}}\bigr).
\end{equation}
Furthermore, if
\[
{nm\lambda^{1/(2r)}}\to\infty,
\]
then
%
%
\begin{equation}
\label{eqerr3}
\Vert\hat{g}-\tilde{g}\Vert^2_{\Lcal_2}=o_p
\bigl(n^{-1}+(nm)^{-1}\lambda^{-{1/(2r)}}\bigr).
\end{equation}
Therefore,
\[
\|\hat{g}-g_0\|_{\Lcal_2}^2=O_p\bigl(\lambda
+n^{-1}+(nm)^{-1}\lambda^{-{1/(2r)}}\bigr).
\]
Taking
\[
\lambda\asymp(nm)^{-{2r/(2r+1)}}
\]
yields
\[
\|\hat{g}-g_0\|_{\Lcal_2}^2=O_p\bigl(n^{-1}+(nm)^{-{2r/
(2r+1)}}\bigr).
\]

We now set to establish bounds (\ref{eqerr1})--(\ref{eqerr3}). For
brevity, we shall assume in what follows that all expectations are
taken conditionally on $m_1,\ldots,m_n$ unless otherwise indicated. Define
%
%
\begin{equation}
\|g\|_{\alpha}^2=\sum_{\nu\ge1} (1+\rho_\nu^{-1}
)^\alpha g_{\nu}^2,
\end{equation}
where $0\le\alpha\le1$.

We begin with $\bar{g}-g_0$.
Write
%
%
\begin{equation}
g_0(\cdot)=\sum_{k\ge1} a_{k}\phi_{k}(\cdot),\qquad g(\cdot)=\sum
_{k\ge1} b_{k}\phi_{k}(\cdot).
\end{equation}
Then
%
%
\begin{equation}
\ell_\infty(g)=\bE\bigl([Y_{11}-g_0(T_{11})]^2
\bigr)+\sum_{k\ge1}(b_{k}-{a}_{k})^2.
\end{equation}
It is not hard to see
%
%
\begin{equation}
\bar{b}_k:=\langle\bar{g},\phi_k\rangle_{\Lcal_2}=\argmin\{
(b_k-a_k)^2+\lambda\rho_k^{-1}b_k^2\}={a_k\over1+\lambda\rho
_k^{-1}}.
\end{equation}
Hence,
\begin{eqnarray*}
\|\bar{g}-g_0\|^2_{\alpha}&=&\sum_{k\ge1} (1+\rho
_{k}^{-1})^{\alpha}(\bar{b}_{k}-a_{k})^2\\
&=&\sum_{k\ge1} (1+\rho_{k}^{-1})^{\alpha}\biggl({ \lambda\rho
_{k}^{-1}\over1+\lambda\rho_{k}^{-1}}\biggr)^2a_{k}^2\\
&\le&c_0\lambda^2\sup_{k\ge1} {\rho_{k}^{-(1+\alpha)}\over
(1+\lambda\rho_{k}^{-1})^2}\sum_{k=1}^\infty\rho
_{k}^{-1}a_{k}^2\\
&\le&c_0\lambda^{1-\alpha} \int\bigl[g_0^{(r)}\bigr]^2.
\end{eqnarray*}

Next, we consider $\tilde{g}-\bar{g}$. Notice that $D\ell
_{mn,\lambda}(\bar{g})=D\ell_{mn,\lambda}(\bar{g})-D\ell_{\infty
,\lambda}(\bar{g})=D\ell_{mn}(\bar{g})-D\ell_{\infty}(\bar{g})$.
Therefore
\begin{eqnarray*}
\bE[D\ell_{mn,\lambda}(\bar{g})f]^2&=&\bE[D\ell
_{mn}(\bar{g})f-D\ell_{\infty}(\bar{g})f]^2\\
&=&{4\over n^2}\sum_{i=1}^n{1\over m_i^2}\var\Biggl[\sum
_{j=1}^{m_i}\bigl([Y_{ij}-\bar{g}(T_{ij})
]f(T_{ij})\bigr)\Biggr].
\end{eqnarray*}
Note that
\begin{eqnarray*}
&&\var\Biggl[\sum_{j=1}^{m_i}\bigl([Y_{ij}-\bar
{g}(T_{ij})]f(T_{ij})\bigr)\Biggr]\\
&&\qquad=\var\Biggl[\bE\Biggl(\sum_{j=1}^{m_i}[Y_{ij}-\bar
{g}(T_{ij})]f(T_{ij})| T\Biggr)\Biggr]+\bE\Biggl[\var
\Biggl(\sum_{j=1}^{m_i}Y_{ij}f(T_{ij})|T\Biggr)\Biggr]\\
&&\qquad=\var\Biggl[\sum_{j=1}^{m_i}\bigl([g_0(T_{ij})-\bar
{g}(T_{ij})]f(T_{ij})\bigr)\Biggr]+\bE\Biggl[\var
\Biggl(\sum_{j=1}^{m_i}Y_{ij}f(T_{ij})|m_i,T\Biggr)\Biggr].
\end{eqnarray*}

The first term on the rightmost-hand side can be bounded by
\begin{eqnarray*}
&&\var\Biggl[\sum_{j=1}^{m_i}\bigl([g_0(T_{ij})-\bar
{g}(T_{ij})]f(T_{ij})\bigr)\Biggr]\\
&&\qquad= m_i\var\bigl([g_0(T_{i1})-\bar{g}(T_{i1})
]f(T_{i1})\bigr)\\
&&\qquad\le m_i\bE\bigl([g_0(T_{i1})-\bar{g}(T_{i1})
]f(T_{i1})\bigr)^2\\
&&\qquad= m_i\int_\Tcal\bigl([g_0(t)-\bar{g}(t)]f(t)
\bigr)^2\,dt\\
&&\qquad\le m_i\int_\Tcal[g_0(t)-\bar{g}(t)]^2\,dt\int_\Tcal
f^2(t)\,dt,
\end{eqnarray*}
where the last inequality follows from the Cauchy--Schwarz inequality.
Together with (\ref{eqerr1}), we get
%
%
\begin{equation}
\var\Biggl[\sum_{j=1}^{m_i}\bigl([g_0(T_{ij})-\bar
{g}(T_{ij})]f(T_{ij})\bigr)\Biggr]\le c_0m_i\|f\|_{\Lcal
_2}^2\lambda.
\end{equation}

We now set out to compute the the second term. First observe that
\begin{eqnarray*}
&&\var\Biggl(\sum_{j=1}^{m_i}Y_{ij}f(T_{ij})|T\Biggr)\\
&&\qquad=\sum
_{j,k=1}^{m_i} f(T_{ij})f(T_{ik})\bigl(C_0(T_{ij}, T_{ik})+\sigma
_0^2\delta_{jk}\bigr),
\end{eqnarray*}
where $\delta_{jk}$ is Kronecker's delta. Therefore,
\begin{eqnarray*}
&&\bE\Biggl[\var\Biggl(\sum_{j=1}^mY_{1j}f(T_{1j})|T
\Biggr)\Big| m_i\Biggr]
\\
&&\qquad=m_i(m_i-1)\int_{\Tcal\times\Tcal} f(s)C_0(s,t)f(t)\,ds\,dt\\
&&\qquad\quad{}+m_i\sigma
_0^2 \|f\|_{\Lcal_2}^2
+m_i\int_{\Tcal} f^2(s)C(s,s)\,ds.
\end{eqnarray*}

Summing up, we have
%
%
\begin{equation}
\bE[D\ell_{mn,\lambda}(\bar{g})\phi_{k}]^2 \le
{c_0\over n^2}\Biggl(\sum_{i=1}^n {1\over m_i}\Biggr)+{4c_k\over n},
\end{equation}
where
%
%
\begin{equation}
c_k=\int_{\Tcal\times\Tcal} \phi_k(s)C_0(s,t)\phi_k(t)\,ds\,dt.
\end{equation}
Therefore,
\begin{eqnarray*}
\bE\Vert\tilde{g}-\bar{g}\Vert_\alpha^2&=&\bE
\biggl\Vert{1\over2}G_\lambda^{-1}D\ell_{nm,\lambda}(\bar{g})
\biggr\Vert_\alpha^2\\
&=& {1\over4}\bE\biggl[\sum_{k\ge1}(1+\rho_{k}^{-1})^\alpha
(1+\lambda\rho_{k}^{-1})^{-2}(D\ell_{nm,\lambda}(\bar
{g})\phi_{k})^2\biggr]\\
&\le& {c_0\over n^2}\Biggl(\sum_{i=1}^n {1\over m_i}\Biggr)\sum
_{k\ge1}(1+\rho_{k}^{-1})^\alpha(1+\lambda\rho_{k}^{-1})^{-2}\\
&&{} +{1\over n}\sum_{k\ge1}(1+\rho_{k}^{-1})^\alpha(1+\lambda
\rho_{k}^{-1})^{-2}c_k.
\end{eqnarray*}

Observe that
\[
\sum_{k\ge1}(1+\rho_{k}^{-1})^\alpha(1+\lambda\rho
_{k}^{-1})^{-2}\le c_0\lambda^{-\alpha-1/(2r)}
\]
and
\[
\sum_{k\ge1}(1+\rho_{k}^{-1})^\alpha(1+\lambda\rho
_{k}^{-1})^{-2}c_k \le\sum_{k\ge1} (1+\rho_{k}^{-1})c_{k}=\bE\|
X\|_{\Wcal_2^r}^2<\infty.
\]
Thus,
\[
\bE\Vert\tilde{g}-\bar{g}\Vert_\alpha^2
\le c_0\Biggl[{1\over n^2}\Biggl(\sum_{i=1}^n {1\over m_i}
\Biggr)\lambda^{-\alpha-1/(2r)}+{1\over n}\Biggr].
\]

It remains to bound $\hat{g}-\tilde{g}$. It can be easily verified that
%
%
\begin{equation}
\hat{g}-\tilde{g}=\tfrac{1}{2}G_\lambda^{-1}[D^2\ell_{\infty
}(\bar{g})(\hat{g}-\bar{g})-D^2\ell_{mn}(\bar{g})(\hat{g}-\bar
{g})].
\end{equation}
Then
\begin{eqnarray*}
\Vert\hat{g}-\tilde{g}\Vert^2_\alpha&=&\sum_{k\ge1}
(1+\rho_{k}^{-1})^{\alpha}(1+\lambda\rho_{k}^{-1})^{-2}\\
&&\hspace*{12.5pt}{} \times\Biggl[{1\over n}\sum_{i=1}^n{1\over m_i}\sum
_{j=1}^{m_i}\bigl(\hat{g}(T_{ij})-\bar{g}(T_{ij})\bigr)
\phi_{k}(T_{ij})\\
&&\hspace*{76.1pt}{}- \int
_{\Tcal}\bigl(\hat{g}(s)-\bar{g}(s)\bigr)\phi_{k}(s)\,ds\Biggr]^2.
\end{eqnarray*}
Clearly, $(\hat{g}-\bar{g})\phi_{k}\in\Hcal(K)$. Write
%
%
\begin{equation}
(\hat{g}-\bar{g})\phi_{k}=\sum_{j\ge1} h_{j}\phi_{j}.
\end{equation}
Then, by the Cauchy--Schwarz inequality,
\begin{eqnarray*}
&&\Biggl[ {1\over n}\sum_{i=1}^n{1\over m_i}\sum_{j=1}^{m_i}\bigl(\hat
{g}(T_{ij})-\bar{g}(T_{ij})\bigr)\phi_{k}(T_{ij})- \int_{\Tcal}\bigl(\hat
{g}(s)-\bar{g}(s)\bigr)\phi_{k}(s)\,ds\Biggr]^2\\
&&\qquad=\Biggl[\sum_{k_1\ge1} h_{k_1}\Biggl({1\over n}\sum
_{i=1}^n{1\over m_i}\sum_{j=1}^{m_i}\phi_{k_1}(T_{ij})- \int_{\Tcal
}\phi_{k_1}(s)\,ds\Biggr)\Biggr]^2\\
&&\qquad\le\biggl[\sum_{k_1\ge1} (1+\rho_{k_1}^{-1})^{\gamma
}h^2_{k_1}\biggr]\\
&&\qquad\quad{}\times\Biggl[\sum_{k_1\ge1} (1+\rho_{k_1}^{-1})^{-\gamma
}\Biggl({1\over n}\sum_{i=1}^n{1\over m_i}\sum_{j=1}^{m_i}\phi
_{k_1}(T_{ij})- \int_{\Tcal}\phi_{k_1}(s)\,ds\Biggr)^2\Biggr]\\
&&\qquad\le\|\hat{g}-\bar{g}\|_\gamma^2(1+\rho_k^{-1})^\gamma\\
&&\qquad\quad{}\times
\Biggl[\sum_{k_1\ge1} (1+\rho_{k_1}^{-1})^{-\gamma}\Biggl({1\over n}\sum
_{i=1}^n{1\over m_i}\sum_{j=1}^{m_i}\phi_{k_1}(T_{ij})- \int_{\Tcal
}\phi_{k_1}(s)\,ds\Biggr)^2\Biggr]
\end{eqnarray*}
for any $0\le\gamma\le1$, where in the last inequality, we used the
fact that
%
%
\begin{equation}
\|(\hat{g}-\bar{g})\phi_{k}\|_\gamma\le\|\hat{g}-\bar{g}\|
_\gamma\|\phi_{k}\|_\gamma=\|\hat{g}-\bar{g}\|_\gamma(1+\rho
_k^{-1})^{\gamma/2}.
\end{equation}
Following a similar calculation as before, it can be shown that
\begin{eqnarray*}
&&\bE\Biggl[\sum_{k_1\ge1} (1+\rho_{k_1}^{-1})^{-\gamma}
\Biggl({1\over n}\sum_{i=1}^n{1\over m_i}\sum_{j=1}^{m_i}\phi
_{k_1}(T_{ij})- \int_{\Tcal}\phi_{k_1}(s)\,ds\Biggr)^2\Biggr]\\
&&\qquad\le{1\over n^2}\Biggl(\sum_{i=1}^n {1\over m_i}\Biggr)\sum_{k_1\ge
1} (1+\rho_{k_1}^{-1})^{-\gamma},
\end{eqnarray*}
which is finite whenever $\gamma>1/2r$. Recall that $m$ is the
harmonic mean of $m_1,\ldots,m_n$. Therefore,
%
%
\begin{equation}
\label{eqgaal}
\Vert\hat{g}-\tilde{g}\Vert^2_\alpha\le O_p
\biggl({1\over nm\lambda^{\alpha+\gamma+1/(2r)}}\biggr)\|\hat{g}-\bar{g}\|
_\gamma^2.
\end{equation}

If $\alpha>1/2r$, then taking $\gamma=\alpha$ yields
%
%
\begin{equation}
\Vert\hat{g}-\tilde{g}\Vert^2_\alpha=O_p\biggl({1\over
nm\lambda^{2\alpha+1/(2r)}}\biggr)\Vert\hat{g}-\bar{g}
\Vert^2_\alpha=o_p(\Vert\hat{g}-\bar{g}\Vert
_\alpha),
\end{equation}
assuming that
%
%
\begin{equation}
nm\lambda^{2\alpha+1/(2r)}\to\infty.
\end{equation}
Together with the triangular inequality
%
%
\begin{equation}
\Vert\tilde{g}-\bar{g}\Vert_\alpha\ge\Vert
\hat{g}-\bar{g}\Vert_\alpha-\Vert\hat{g}-\tilde
{g}\Vert_\alpha=\bigl(1-o_p(1)\bigr)\Vert\hat{g}-\bar{g}
\Vert_\alpha.
\end{equation}
Therefore,
%
%
\begin{equation}
\Vert\hat{g}-\bar{g}\Vert_\alpha^2=O_p(
\Vert\tilde{g}-\bar{g}\Vert_\alpha^2)=O_p
\bigl(n^{-1}+(nm)^{-1}\lambda^{-\alpha-{1/(2r)}}\bigr).
\end{equation}
Together with (\ref{eqgaal}),
\begin{eqnarray*}
\Vert\hat{g}-\tilde{g}\Vert^2_{\Lcal_2}&=&O_p
\biggl({1\over nm\lambda^{\alpha+1/(2r)}}\bigl(n^{-1}+(nm)^{-1}\lambda
^{-\alpha-{1/(2r)}}\bigr)\biggr)\\
&=&o_p\bigl(n^{-1}\lambda^{\alpha}+(nm)^{-1}\lambda^{-{1/
(2r)}}\bigr).
\end{eqnarray*}

We conclude by noting that in the case when $m_1,\ldots, m_n$ are
random, $m$ can also be replaced with the expectation of the harmonic
mean thanks to the law of large numbers.
\end{pf*}
%
%
%
%
\begin{pf*}{Proof of Proposition \ref{prequiv}}
Let $Q_{T,\lambda}$ be the smoothing spline operator, that is,
$Q_{T,\lambda}(f_1,\ldots, f_m)$ is the solution to
\[
\min_{g\in\Wcal_2^r}\Biggl\{{1\over m}\sum_{j=1}^m
\bigl(f_j-g(T_{j})\bigr)^2 +\lambda\int_\Tcal\bigl[g^{(r)}(t)
\bigr]^2\,dt\Biggr\}.
\]
It is clear that
\[
\hat{g}_\lambda=Q_{T,\lambda}(\bar{Y}_{\cdot1}, \bar{Y}_{\cdot
2},\ldots, \bar{Y}_{\cdot m})
\]
and
\[
\tilde{X}_{i,\lambda}=Q_{T,\lambda}(Y_{i1}, {Y}_{i2},\ldots, {Y}_{im}).
\]
Because $Q_{T,\lambda}$ is a linear operator [see, e.g., \citet{w90}],
we have
\begin{eqnarray*}
\tilde{g}_\lambda&=&{1\over n}\sum_{i=1} \tilde{X}_{i,\lambda
}={1\over n}\sum_{i=1}^nQ_{T,\lambda}(Y_{i1}, {Y}_{i2},\ldots,
{Y}_{im})\\
&=&Q_{T,\lambda}(\bar{Y}_{\cdot1}, \bar{Y}_{\cdot2},\ldots
, \bar{Y}_{\cdot m})=\hat{g}_\lambda.
\end{eqnarray*}
\upqed\end{pf*}

\section*{Acknowledgments}

We thank an Associate Editor and two referees for their constructive
comments which have helped to improve the presentation of the paper.


%

\printaddresses

\end{document}